\documentclass[11pt,a4]{article}
\usepackage{amsfonts}
\usepackage{graphics,amssymb,amsthm,amsmath,epsf}
\usepackage{graphicx}
\usepackage{subfigure}

\newtheorem{algorithm}{Algorithm}[section]
\newtheorem{theorem}{Theorem}[section]

\numberwithin{equation}{section}
\newtheorem{remark}{Remark}[section]

\usepackage{fancyhdr}

\begin{document}

\vspace{5cm}

\title{Generalized inference for the common mean of several lognormal populations}
\author{J. Behboodian* and A. A. Jafari**
\\ *Department of Mathematics, Shiraz Islamic Azad University \\
Shiraz, IRAN \\
email: Behboodian@stat.susc.ac.ir\\
*Department of Statistics, Shiraz University, Shiraz, IRAN}
\date{}
\maketitle

\begin{abstract}
A hypothesis testing and an interval estimation are studied for the common
mean of several lognormal populations. Two methods are given based on the
concept of generalized p-value and generalized confidence interval. These new
methods are exact and can be used without restriction on sample sizes, number
of populations, or difference hypotheses. A simulation study for coverage
probability, size and power shown that the new methods are better than the
existing methods. A numerical example is given with some real medical data.
\end{abstract}

\bigskip

KewWords: Lognormal population, Common mean, Generalized variable,
Generalized p-value, Generalized confidence interval.

\section{Introduction}

The statistical analysis that combines the results of several independent is known as meta-analysis and it is used in clinical
trails and behavioral sciences.

Consider  we have $k$ independent normal populations with means
$a\mu +b\sigma _{i}^{2}$ and variances $\sigma _{i}^{2}$. Also we
have a random samples of sizes $n_{i}$, $i=1,...,k$ from each one.
We denote these samples by $Y_{ij}\sim N(a\mu +b\sigma
_{i}^{2},\sigma _{i}^{2}),$ $i=1,...,k,$ $ j=1,...,n_{i}$, where
$a\neq 0,$ and $b$ are constant. The problem of interest is to
combine the summary statistics of samples for statistical
inference about the parameter $\mu $. The statistical analysis
that combines the results of several independent used in clinical
trails and behavioral sciences.

If $a=1$ and $b=0$ then, $Y_{ij}\sim N(\mu ,\sigma _{i}^{2})$ and
this problem is known as the common mean for several normal
populations. There are some inference for this problem in
statistical literature. For example see; Krishnamoorthy and Lu
(2003), Lin and Lee (2005). If $a=1$ and $b=-0.5,$ then
$Y_{ij}\thicksim N(\mu -0.5\sigma _{i}^{2},\sigma _{i}^{2})$ and
this is equivalent to problem of common mean of several lognormal
populations. Our interest in this paper is inference about this
problem. For the common lognormal mean, a few authors proposed
approximate methods: Ahmed et al (2001) proposed an estimator and
approximate confidence interval for the common lognormal mean;
Baklizi and Ebrahem (2005) studied several types of large samples
and bootstrap intervals; Gupta and Li (2005) developed procedures
for estimating the common mean and investigated the performance of
the resulting confidence interval for two lognormal populations.

In this paper, we first propose estimation of $\mu$ when the
variances, $\sigma_i^2$ are known. Then two  methods are given
that are applicable for both hypothesis testing and interval
estimation for $\mu$, based on the concepts of generalized
$p$-value and generalized confidence interval. These methods are
based on extending the method of Krishnamoorthy and Lu (2003) and
the method of Lin and Lee (2005), which are used for the problem
of common mean of several normal populations. Our methods also are
applicable for the common mean of several lognormal for the
interval mean of $k$ lognormal populations. This cahpter also is
devoted to a short review regarding the existing method for
inference of the common lognormal mean and application of our two
methods for this problem. Finally, we give a numerical example for
the common lognormal mean and by Monte Carlo simulation, we
compare the coverage probabilities, size and power of these
methods for the common mean of two lognormal populations.

\begin{theorem}
Let $Y_{ij}\thicksim N(a\mu +b\sigma _{i}^{2},\sigma _{i}^{2})$,
$i=1,...,k$, $j=1,...,n_{i}$, where $a\neq 0,$ $b$ are constants
and $\sigma _{i}^{2}$'s are known. The estimator
\begin{equation}\label{eq2.1.1}
\hat{\mu}=\dfrac{\sum\limits_{i=1}^{k}\dfrac{n_{i}\bar{Y}_{i.}}{\sigma
_{i}^{2}}-nb}{a\sum\limits_{i=1}^{k}\dfrac{n_{i}}{\sigma
_{i}^{2}}}
\end{equation}
is UMVUE and MLE for $\mu $ and $\hat{\mu}\thicksim N(\mu
,1/(a^{2}\sum\limits_{i=1}^{k}(\dfrac{n_{i}}{\sigma _{i}^{2}}))).$
\end{theorem}

\noindent\begin{proof}
 The probability density function for $Y_{ij} $ is
\[
f_{Y_{ij}}(y_{ij})=(2\pi \sigma
_{i}^{2})^{\frac{-1}{2}}e^{\frac{-1}{2\sigma _{i}^{2}}a^{2}\mu
^{2}}\times e^{\frac{-1}{2\sigma _{i}^{2}}(y_{ij}-b\sigma
_{i}^{2})^{2}}\times e^{\frac{a\mu }{\sigma
_{i}^{2}}(y_{ij}-b\sigma _{i}^{2})}.
\]

Since the distribution of $Y_{ij}$ is from exponential family, in
the form $A(\mu )B(y)e^{C(\mu )D(y)}$, then
\[
T=\sum\limits_{i=1}^{k}\sum\limits_{j=1}^{n_{i}}\frac{1}{\sigma _{i}^{2}}%
(Y_{ij}-b\sigma _{i}^{2})=\sum\limits_{i=1}^{k}\dfrac{n_{i}\bar{Y}_{i.}}{%
\sigma _{i}^{2}}-nb
\]
is UMVUE\ for $E(T)=a\mu
\sum\limits_{i=1}^{k}\dfrac{n_{i}}{\sigma _{i}^{2}}$ and $\
\hat{\mu}=T/\sum\limits_{i=1}^{k}\dfrac{an_{i}}{\sigma _{i}^{2}}$
is UMVUE for $\mu $ (see Casella and Berger, 1990, page 263)$.$
It is easy to prove the rest of the theorem.
\end{proof}

\begin{remark}
If b=0 then $\hat{\mu}$ is the best linear unbiased estimator for
$\mu .$
\end{remark}

\begin{remark} If $Y_{ij}=\ln (X_{ij})\thicksim N(\mu
-0.5\sigma
_{i}^{2},\sigma _{i}^{2}),$ i.e. $X_{ij}$ is a lognormal variable, then $%
T=\exp (\hat{\mu}-1/\sum\limits_{i=1}^{k}\dfrac{2n_{i}}{\sigma
_{i}^{2}})$
is UMVUE for $E(X_{ij})=e^{\mu }$, but the MLE of $e^{\mu }$ is $e^{\hat{\mu}%
}.$
\end{remark}

\begin{remark} If $\sigma _{i}^{2}$ are unknown, then
we cannot find a closed form for MLE's of $\mu $; we have to use a
numerical approximation.
\end{remark}

\section{ Generalized inferences for $\mu $}
Suppose $Y_{ij}\thicksim N(a\mu +b\sigma _{i}^{2},\sigma _{i}^{2})$%
, $i=1,...,k$, $j=1,...,n_{i}$, where $a\neq 0,$ $b$ are constants. For the $%
i$th population, let
\[
\bar{Y}_{i.}=\frac{1}{n_{i}}\sum\limits_{i=1}^{n_{i}}Y_{ij}\text{
\ \ \ \ ,
\ \ \ }S_{i}^{2}=\frac{1}{n_{i}-1}\sum\limits_{i=1}^{n_{i}}(Y_{ij}-\bar{Y}%
_{i.})^{2},
\]
be the sample mean and sample variance.

In this section, by using the idea of generalized $p$-value and by
extending (i) the method of Krishnamoorthy and Lu (2003) and (ii)
the \ method of \ Lin and Lee (2005), for the problem of common
mean of normal populations, we give two generalized pivot
variables for interval estimation and hypothesis testing for $\mu
$ and we obtain two generalized $p$-values for testing hypothesis
\begin{equation}\label{eq2.2.1}
H_{\circ }:\mu \leq \mu _{\circ }\text{ \ \ \ vs \ \ \ \
}H_{1}:\mu <\mu _{\circ }.
\end{equation}

\subsection{A weighted linear combination}

It is clear that $\bar{Y}_{i.}\thicksim N(a\mu +b\sigma
_{i}^{2},\sigma _{i}^{2}/n_{i})$, $i=1,...,k$. Therefore, the
generalized pivot variable for estimating $\mu $ based on the
$i$th sample is
\begin{eqnarray}\label{eq2.2.2}
T_{i}^{\ast } &=&\frac{1}{a}(\bar{y}_{i.}-b\frac{(n_{i}-1)s_{i}^{2}}{U_{i}}%
-Z_{i}\sqrt{\frac{(n_{i}-1)s_{i}^{2}}{n_{i}U_{i}}}) \\
&=&\frac{1}{a}(\bar{y}_{i.}-b\frac{s_{i}^{2}}{S_{i}^{2}}\sigma _{i}^{2}-Z_{i}%
\sqrt{\frac{s_{i}^{2}}{n_{i}S_{i}^{2}}\sigma _{i}^{2}})\text{,}
\nonumber
\end{eqnarray}
where
$$Z_{i}=\dfrac{\bar{Y}_{i.}-(a\mu +b\sigma _{i}^{2})}{\sqrt{\sigma
_{i}^{2}/n_{i}}}\sim N(0,1), \ \ \
U_{i}=\dfrac{(n_{i}-1)S_{i}^{2}}{\sigma _{i}^{2}}\sim \chi
_{(n_{i}-1)}^{2},$$ and $(\bar{y}_{i.}$,$s_{i}^{2})$ is the
observed value of $(\bar{Y}_{i.}$, $S_{i}^{2})$.

The generalized pivot variable for estimating $\sigma _{i}^{2}$
based on the $i$th sample is given by
\begin{equation}\label{eq2.2.3}
R_{i}=\frac{(n_{i}-1)s_{i}^{2}}{V_{i}}=\frac{s_{i}^{2}}{S_{i}^{2}}\sigma
_{i}^{2}\text{, \ }i=1,...,k\text{,}
\end{equation}
where $V_{i}=(n_{i}-1)S_{i}^{2}/\sigma _{i}^{2}$ are independent
$\chi _{(n_{i}-1)}^{2}$ random variables (Weerahandi, 1995).

The generalized variable that we want to propose is a weighted
average of the generalized pivot variables $T_{i}^{\ast }$ in
(\ref{eq2.2.2}). The weights are inversely proportional to the
generalized pivot variables $R_{i}$ in (\ref{eq2.2.3}) for the
variances, and they are directly proportional to the sample
sizes. (see Krishnamoorthy and Lu, 2003).

Let $\bar{Y}=(\bar{Y}_{1.},...,\bar{Y}_{k.})$ and
$V=(V_{1},...,V_{k})$, with the observed values  $\bar{y}$ and
$v$, respectively. Then, the generalized variable can be
expressed as
\begin{eqnarray}\label{eq2.2.4}
T(\bar{Y},V;\bar{y},v) &=&\frac{\sum\limits_{i=1}^{k}\dfrac{n_{i}V_{i}}{%
(n_{i}-1)s_{i}^{2}}\left[ \bar{y}_{i.}-b\dfrac{(n_{i}-1)s_{i}^{2}}{U_{i}}%
-Z_{i}\sqrt{\dfrac{(n_{i}-1)s_{i}^{2}}{n_{i}U_{i}}}\right] }{%
a\sum\limits_{j=1}^{k}\dfrac{n_{j}V_{j}}{(n_{j}-1)s_{j}^{2}}}-\mu  \\
&=&\sum\limits_{i=1}^{k}W_{i}T_{i}^{\ast }-\mu ,  \nonumber
\end{eqnarray}
where the weights are
\[
W_{i}=\frac{\dfrac{n_{i}V_{i}}{(n_{i}-1)s_{i}^{2}}}{\sum\limits_{j=1}^{k}%
\dfrac{n_{j}V_{j}}{(n_{j}-1)s_{j}^{2}}}\text{, \ }i=1,...,k.
\]

The distribution of $T(\bar{Y},V;\bar{y},v)$ is an increasing
function with respect to $\mu $. Therefore, the generalized
$p$-value for (\ref{eq2.2.2}) is given by
\begin{eqnarray}\label{eq2.2.5}
p &=&P(T(\bar{Y},V;\bar{y},v)\leqslant T(\bar{y},v;\bar{y},v)\mid
\mu
=\mu _{\circ }) \\
&=&P(\sum\limits_{i=1}^{k}W_{i}T_{i}^{\ast }\leqslant \mu _{\circ
}). \nonumber
\end{eqnarray}

This generalized $p$-value can be well approximated by a Monte
Carlo simulation using the following algorithm:

\begin{algorithm} \label{alg2.1}
 For a given $(n_{1},...,n_{k})$, $\bar{Y}=(\bar{y}%
_{1.},...,\bar{y}_{k.})$ and ($s_{1}^{2},...,s_{k}^{2})$:

\noindent For $j=1,m$

\noindent Generate $U_{l}\thicksim \chi _{(n_{l}-1)}^{2},$ \
$l=1,...,k$

\noindent Generate $V_{l}\thicksim \chi _{(n_{l}-1)}^{2},$ \
$l=1,...,k$

\noindent Generate $Z_{l}\thicksim N(0,1),$ \ $l=1,...,k$

\noindent Compute $W_{1},...,W_{k}$

\noindent Compute $T_{j}=\sum\limits_{l=1}^{k}W_{l}T_{l}^{\ast }$

\noindent (end $j$ loop)

Let $\gamma_{j}=1$ if $T_{j}\leqslant \mu _{\circ }$, else
$k_{j}=0$. Then $\dfrac{1}{m}\sum\limits_{j=1}^{m}\gamma_{j}$ is a
Monte Carlo estimate of the generalized $p$-value for
(\ref{eq2.2.5}).
\end{algorithm}

\begin{remark} $T^{\ast }=\sum\limits_{i=1}^{k}W_{i}T_{i}^{\ast
}$ is a generalized pivot variable for $\mu $ and we can use that
to obtain a generalized confidence interval for $\mu$.
\end{remark}

\begin{remark}
If $a=1$ and $b=0,$ then
\begin{equation}
T(\bar{Y},V;\bar{y},v)=\frac{\sum\limits_{i=1}^{k}\dfrac{n_{i}V_{i}}{%
(n_{i}-1)s_{i}^{2}}\left[ \bar{y}_{i.}-Z_{i}\sqrt{\dfrac{(n_{i}-1)s_{i}^{2}}{%
n_{i}U_{i}}}\right] }{\sum\limits_{j=1}^{k}\dfrac{n_{j}V_{j}}{%
(n_{j}-1)s_{j}^{2}}}-\mu
\end{equation}
and this generalized variable is introduced by Krishnamoorthy and
Lu (2003) for inference on the common mean of several normal
populations.
\end{remark}

\subsection{A generalized variable based on UMVUE}

From theorem 1, we have
\[
Z=\mid a\mid \sqrt{\sum\limits_{i=1}^{k}\frac{n_{i}}{\sigma _{i}^{2}}}(\hat{%
\mu}-\mu )\thicksim N(0,1).
\]

We know that $R_{i}=\dfrac{(n_{i}-1)s_{i}^{2}}{U_{i}}$ is a
generalized pivot variable for $\sigma _{i}^{2},$ $i=1,...,k,$
where $U_{i}\thicksim \chi _{(n_{i}-1)}^{2}.$

Let $\bar{Y}=(\bar{Y}_{1.},...,\bar{Y}_{k.})$ and
$U=(U_{1},...,U_{k})$, with the observed values $\bar{y}$ and
$u$, respectively. We define a generalized variable for $\mu $
based on the UMVUE for $\mu $ in (\ref{eq2.1.1}) by
\begin{eqnarray}\label{eq2.2.7}
T(\bar{Y},U;\bar{y},u) &=&\frac{\sum\limits_{i=1}^{k}\dfrac{n_{i}\bar{y}_{i}%
}{(n_{i}-1)s_{i}^{2}}U_{i}-nb}{a\sum\limits_{j=1}^{k}\dfrac{n_{j}}{%
(n_{j}-1)s_{j}^{2}}U_{j}}-\frac{Z}{\mid a\mid \sqrt{\sum\limits_{j=1}^{k}%
\dfrac{n_{j}}{(n_{j}-1)s_{j}^{2}}U_{j}}}-\mu  \\
&=&\frac{\sum\limits_{i=1}^{k}\dfrac{n_{i}\bar{y}_{i}}{\sigma _{i}^{2}}%
\dfrac{S_{i}^{2}}{s_{i}^{2}}-nb}{a\sum\limits_{j=1}^{k}\dfrac{n_{j}}{\sigma
_{j}^{2}}\dfrac{S_{j}^{2}}{s_{j}^{2}}}-\frac{\sqrt{\sum\limits_{i=1}^{k}%
\dfrac{n_{i}}{\sigma _{i}^{2}}}(\hat{\mu}-\mu )}{\sqrt{\sum\limits_{j=1}^{k}%
\dfrac{n_{j}}{\sigma _{j}^{2}}\dfrac{S_{j}^{2}}{s_{j}^{2}}}}-\mu
.  \nonumber
\end{eqnarray}

The distribution of $T(\bar{Y},U;\bar{y},u)$ is an increasing
function with respect to $\mu $, and therefore the generalized
$p$-value for testing (\ref{eq2.2.1}) is
\begin{eqnarray}
p&=&P(T(\bar{Y},U;\bar{y},u)\leqslant T(\bar{y},u;\bar{y},u)\mid
\mu
=\mu _{\circ })=P(T^{\ast }\leqslant \mu _{\circ }) \\
&=&1-E\left[ \Phi (\frac{\mid a\mid }{a}\frac{\sum\limits_{i=1}^{k}\dfrac{%
n_{i}\bar{y}_{i}}{(n_{i}-1)s_{i}^{2}}U_{i}-nb}{\sqrt{\sum\limits_{j=1}^{k}%
\dfrac{n_{j}}{(n_{j}-1)s_{j}^{2}}U_{j}}}-\mid a\mid \sqrt{%
\sum\limits_{j=1}^{k}\dfrac{n_{j}}{(n_{j}-1)s_{j}^{2}}U_{j}}\mu _{\circ })%
\right] ,  \nonumber
\end{eqnarray}
where
\begin{equation}
T^{\ast }=\frac{\sum\limits_{i=1}^{k}\dfrac{n_{i}\bar{y}_{i}}{\sigma _{i}^{2}%
}\dfrac{S_{i}^{2}}{s_{i}^{2}}-bn}{a\sum\limits_{j=1}^{k}\dfrac{n_{j}}{\sigma
_{j}^{2}}\dfrac{S_{j}^{2}}{s_{j}^{2}}}-\frac{\sqrt{\sum\limits_{i=1}^{k}%
\dfrac{n_{i}}{\sigma _{i}^{2}}}(\hat{\mu}-\mu )}{\sqrt{\sum\limits_{j=1}^{k}%
\dfrac{n_{j}}{(n_{j}-1)s_{j}^{2}}U_{j}}},
\end{equation}
and $\Phi $ is distribution function of the standard normal
variable and
expectation is taken with respect to chi-square random variables with $%
n_{i}-1,$ $i=1,...,k,$ degrees of freedom.

This generalized $p$-value can be well approximated by a Monte
Carlo simulation like the algorithm \ref{alg2.1}.

\begin{remark} $T^{\ast }$ in (2.10) is a generalized
pivot variable for $\mu $ and we can use that to obtain a
generalized confidence interval for $\mu$.
\end{remark}

\begin{remark}

If $a=1$ and $b=0,$ then
\[
T(\bar{Y},U;\bar{y},u)=\frac{\sum\limits_{i=1}^{k}\dfrac{n_{i}\bar{y}_{i}}{%
(n_{i}-1)s_{i}^{2}}U_{i}}{\sum\limits_{j=1}^{k}\dfrac{n_{j}}{%
(n_{j}-1)s_{j}^{2}}U_{j}}-\frac{Z}{\sqrt{\sum\limits_{j=1}^{k}\dfrac{n_{j}}{%
(n_{j}-1)s_{j}^{2}}U_{j}}}-\mu ,
\]
which is a generalized variable, introduced by Lin and Lee
(2005), for the common mean of several normal populations.
\end{remark}

\begin{remark} For testing the hypothesis of the form $$ H_{\circ }:\mu =\mu
_{\circ }\text{ \ \ \ vs}\text{ \ \ \ }H_{1}:\mu \neq \mu _{\circ
},$$ the $p$-value is
\begin{equation}
p=2\min \{P\{T^{\ast }<\mu _{\circ }\},P\{T^{\ast }>\mu _{\circ
}\}\},
\end{equation}
and $H_{\circ }$ can be rejected when $p<\alpha $.
\end{remark}

\section{ Methods for Common lognormal mean}

Consider independent  $X_{ij}$ with lognormal distribution, for
$i=1,...,k$, $j=1,...,n_{i}$, and assume that $\theta
_{1}=...=\theta _{k}=\varphi >0$, where $\theta
_{i}=E(X_{ij})=\exp (\mu _{i}+\sigma _{i}^{2})$, i.e., the $k$
lognormal populations have common mean $\varphi .$ Therefore, we have $%
Y_{ij}=\ln (X_{ij})\thicksim N(\mu -0.5\sigma _{i}^{2},\sigma
_{i}^{2})$, where $\mu =\ln \varphi$, and to find a confidence
interval for $\varphi$, it is enough to have a confidence interval
for $\mu$, and  a hypothesis test for $\varphi$ is equivalent to
a hypothesis test for $\mu$. For example the hypothesis test
\[
H_{\circ }:\varphi \leqslant \varphi _{\circ }\text{ \ }vs\text{ \ \ }%
H_{1}:\varphi> \varphi _{\circ },
\]
is equivalent to
\[
H_{\circ }:\mu \leqslant \ln \varphi _{\circ }\text{ \ }vs\text{ \ \ }%
H_{1}:\mu> \ln \varphi _{\circ }.
\]

It is useful to review the existing methods for the problem of
common lognormal mean.

\subsection{ Ahmed method}
 Let $X_{ij}\thicksim LN(\theta ,\tau
_{i}^{2}),$
$i=1,...,m$, $j=1,..,n_{i}.$ Then a combined sample estimate of $%
E(X_{ij})=\theta $ is given by
\[
\tilde{\theta}=\frac{\sum\limits_{i=1}^{m}\dfrac{n_{i}}{\hat{v}_{i}}\hat{%
\theta}_{i}}{\sum\limits_{i=1}^{m}\dfrac{n_{i}}{\hat{v}_{i}}},
\]
where $\hat{v}_{i}=\hat{\sigma}_{i}^{2}(1+0.5\hat{\sigma}_{i}^{2})\exp (2%
\hat{\mu}_{i}+\hat{\sigma}_{i}^{2}),$ $\hat{\theta}_{i}=\exp (\hat{\mu}%
_{i}+0.5\hat{\sigma}_{i}^{2})$, $\hat{\mu}_{i}=\bar{Y}_{i.}$ and $\hat{\sigma%
}_{i}^{2}=\dfrac{n_{i}-1}{n_{i}}S_{i}^{2}.$

The estimator $\tilde{\theta}$ is asymptotically normal with mean
$\theta $ and asymptotic variance
$(\sum\limits_{i=1}^{m}\dfrac{n_{i}}{v_{i}})^{-1}$,
which can be estimated by $(\sum\limits_{i=1}^{m}\dfrac{n_{i}}{\hat{v}_{i}}%
)^{-1}.$ Therefore, a $100(1-\alpha )\%$ confidence interval for
$\theta $ is
\begin{equation}
\tilde{\theta}\pm Z_{\alpha /2}(\sum\limits_{i=1}^{m}\dfrac{n_{i}}{\hat{v}%
_{i}})^{-1/2}.
\end{equation}

\subsection{ Baklizi and Ebrahem method }
 The acceptance set
for all $\theta $ is
\begin{equation}
\sum\limits_{i=1}^{m}\dfrac{n_{i}(\hat{\theta}_{i}-\theta )^{2}}{\hat{v}_{i}}%
\leqslant \chi _{\alpha ,m}^{2}.
\end{equation}

This is a quadratic function in $\theta $ whose two roots can be
found directly. Since the coefficient of $\theta ^{2}$ in this
expression is positive, it follows that the set of all values of
$\theta $ between the two roots is the desired confidence
interval.

\subsection{ Gupta and Li method}
 Let ${\bf \theta }$
$=(\mu ,\sigma _{1},\sigma _{2})$ be a vector of parameters,
where $\mu =\ln \eta =\mu _{i}+0.5\sigma _{i}^{2},$ $i=1,2$ and
$\eta $ is the common mean. The joint log-likelihood function
based on the log-transformed data of two independent log-normal
populations is given by
\begin{eqnarray*}
\ln l({\bf \theta }) &=&(-(n_{1}+n_{2})/2)\ln 2\pi -n_{1}\ln
\sigma
_{1}-n_{2}\ln \sigma _{2}-0.5(t_{1}+t_{2})+\frac{\mu }{\sigma _{1}^{2}}t_{1}-%
\frac{1}{2\sigma _{1}^{2}}t_{3} \\
&&-\frac{(\mu -\sigma _{1}^{2}/2)^{2}n_{1}}{2\sigma _{1}^{2}}+\frac{\mu }{%
\sigma _{2}^{2}}t_{2}-\frac{1}{2\sigma _{2}^{2}}t_{4}-\frac{(\mu
-\sigma _{2}^{2}/2)^{2}n_{2}}{2\sigma _{2}^{2}},
\end{eqnarray*}
where
$$(t_{1},t_{2},t_{3},t_{4})=(\sum\limits_{j}\ln
x_{1j},\sum\limits_{j}\ln x_{2j},\sum\limits_{j}(\ln
x_{1j})^{2},\sum\limits_{j}(\ln x_{2j})^{2}).$$

Let $\hat{\mu}$ is MLE for $\mu .$ The asymptotic variance of
$\hat{\mu}$ is
\[
Var(\hat{\mu})=\frac{(\dfrac{2n_{1}}{\hat{\sigma}_{1}^{2}}+n_{1})(\dfrac{%
2n_{2}}{\hat{\sigma}_{2}^{2}}+n_{2})}{\dfrac{2n_{1}^{2}}{\hat{\sigma}_{1}^{4}%
}(\dfrac{2n_{2}}{\hat{\sigma}_{2}^{2}}+n_{2})+\dfrac{2n_{2}^{2}}{\hat{\sigma}%
_{2}^{4}}(\dfrac{2n_{1}}{\hat{\sigma}_{1}^{2}}+n_{1})},
\]
where $\hat{\sigma}_{1}$ and $\hat{\sigma}_{2}$ are MLEs for
$\sigma _{1}$ and $\sigma _{2}.$ A $100(1-\alpha )\%$ confidence
interval for $\eta =e^{\mu }$ is
\begin{equation}
\exp (\hat{\mu}\pm Z_{\alpha /2}\times SD(\hat{\mu})).
\end{equation}

\subsection{ Generalized inferences} In fact, the problem of
common lognormal mean is a special case of our model when $a=1$
and $b=-\dfrac{1}{2}.$ Thus, the generalized variable in
(\ref{eq2.2.4}) becomes
\[
T(\bar{Y},V;\bar{y},v)=\frac{\sum\limits_{i=1}^{k}\dfrac{n_{i}V_{i}}{%
(n_{i}-1)s_{i}^{2}}\left[ \bar{y}_{i.}+\dfrac{(n_{i}-1)s_{i}^{2}}{2U_{i}}%
-Z_{i}\sqrt{\dfrac{(n_{i}-1)s_{i}^{2}}{n_{i}U_{i}}}\right] }{%
\sum\limits_{j=1}^{k}\dfrac{n_{j}V_{j}}{(n_{j}-1)s_{j}^{2}}}-\mu ,
\]
and the generalized variable in (\ref{eq2.2.7}) becomes
\[
T(\bar{Y},U;\bar{y},u)=\frac{\sum\limits_{i=1}^{k}\dfrac{n_{i}\bar{y}_{i}}{%
(n_{i}-1)s_{i}^{2}}U_{i}+\dfrac{n}{2}}{\sum\limits_{j=1}^{k}\dfrac{n_{j}}{%
(n_{j}-1)s_{j}^{2}}U_{j}}-\frac{Z}{\sqrt{\sum\limits_{j=1}^{k}\dfrac{n_{j}}{%
(n_{j}-1)s_{j}^{2}}U_{j}}}-\mu
\]

\section{Numerical Studies}
In this section, we give a numerical example and compare our
methods with other methods for the problem of common lognormal
mean.
\subsection{An example} The data come from the Regenstrief
Medical Record System (RMRS) (MCDonald et al, 1988; Zhou et al,
1997) on effects of race on medical charges of patients with type
I diabetes who had received inpatient or outpatient care at least
two occasions during the period from 1 January 1993, through 30
June 1994. The data set consists of 119 African American patients
and 106 white patients. The mean medical charges and their
corresponding variance for the African American and white groups
are given in Table \ref{table2.1}.

\begin{table}[ht]
\begin{center}
\caption{Sample means and sample variances of the original and
the log-transformed RMRS data} \label{table2.1}
\begin{tabular}{llll}
\hline Data \ \ \ \ \ \ \ \ \ \ \ \ \ \ \ \ \ \ \  & Patients
group \ \ \ \ \ \ \ \ \  & Sample mean \ \  & Sample variance \
\\ \hline
\begin{tabular}{l}
Original \\
\\
Log-transform
\end{tabular}
&
\begin{tabular}{l}
African American \\
White \\
African American \\
White
\end{tabular}
&
\begin{tabular}{l}
\$18,850 \\
\$18,584 \\
9.06695 \\
8.69306
\end{tabular}
&
\begin{tabular}{l}
26.897$^{\text{2}}$ \\
30.694$^{\text{2}}$ \\
1.824 \\
2.629
\end{tabular}
\\ \hline
\end{tabular}
\end{center}
\end{table}

The studies show that (i) lognormal model adequately describes
the both data sets. (ii) the variances of the two sets are not
equal. (iii) the means of the two sets are equal (see Gupta and
Li, 2005). Therefore, the average medical costs for African
American patients and white patients are the same. We want to
test that this average medical costs is 20,000\$, i.e. the
hypothesis test
\begin{equation}
H_{\circ }:\varphi =20000\text{ \ }vs\text{ \ \ }H_{1}:\varphi
\neq 20000,
\end{equation}

The $p$-values for this test, with different methods are given in
Table \ref{table2.2} and
the confidence intervals are given in Table \ref{table2.3}. Therefore, we cannot reject $%
H_{\circ }.$

\begin{table}[ht]
\begin{center}
\caption{$p$-values for hypothesis test of the common lognormal
mean $ \varphi $} \label{table2.2}
\begin{tabular}{ll}
\hline \ \  Methods \ \ \ \ \ \ \ \ \  & $p$-values \\ \hline
Likelihood Ratio Test & 0.5245 \\
Ahmed  & 0.5582\\
Gupta and Li  & 0.5343\\
 First Generalized $p$-value \ \ \ \ \ \ \ \ \ \ \ \ \ \ \ \ \ \ \ \ \ \
\ \ \ \ \ \ \ \ \ \  & 0.4348\\
 Second Generalized $p$-value & 0.4732 \\ \hline
\end{tabular}
\end{center}
\end{table}

\begin{table}
\begin{center}
\caption{ Interval estimation for the common lognormal mean
$\varphi $} \label{table2.3}
\begin{tabular}{lcc}
\hline \ \ Methods & Intervals & Width \\ \hline
\begin{tabular}{l}
Ahmed  \\
Gupta and Li  \\
Baklizi and Ebrahem  \\
First Generalized confidence \\
Second Generalized confidence
\end{tabular}
&
\begin{tabular}{l}
(15831.21 , 27720.26) \\
(16596.91 , 28658.17) \\
(14372.59 , 29178.79) \\
(17286.30 , 30701.92) \\
(17090.54 , 29998.23)
\end{tabular}
&
\begin{tabular}{l}
11889.14 \\
12061.19 \\
14806.20 \\
13415.62 \\
12907.69
\end{tabular}
\\ \hline
\end{tabular}
\end{center}
\end{table}

\subsection{Simulation study}

A simulation study is performed for inference about  the common
lognormal mean, $\varphi$. The purpose of the simulation is to
compare the size, power and coverage probability of each of the
introduced methods with the  others existing  for two lognormal
populations. For this purpose, several data sets from two normal
distributions, with means $\mu -0.5\sigma _{i}^{2}$ and variances
$\sigma _{i}^{2},$ $i=1,2,$ where $\mu =\ln \varphi , $ were
created. For each condition $10000$ simulations are used. The
sizes are given in table \ref{table2.4}, and the powers in tables
\ref{table2.5} and \ref{table2.6}, and the coverage probability
in tables \ref{table2.7}, \ref{table2.8} and \ref{table2.9}.
These methods are

(1) Likelihood ratio test

(2) Ahmed method

(3) Gupta and Li method

(4) Baklizi and Ebrahem method

(5) First Generalized variable in (\ref{eq2.2.3})

(6) Second Generalized variable in (\ref{eq2.2.7})

 The tables
show that

\begin{itemize}
\item The simulated sizes of the two new methods are satisfactory
since they are close to the significance level, 0.05.
\item The power of the first generalized method  is better than other methods when the
sample sizes are large.
\item The coverage probabilities of our generalized  methods are close to
the significance level and they are better than the coverage
probabilities of existing methods.
\end{itemize}

\begin{table}
\begin{center}
\caption{ Simulated  sizes of the tests for $H_{\circ }:\varphi =1
$\ $vs$ $H_{1}:\varphi \neq 1$ at 5\% significance level when
$\mu =0$ and $\sigma _{1}^{2}=1.$} \label{table2.4}
\begin{tabular}{|c|ccccc|}
\hline
\begin{tabular}{r}
$\sigma _{2}^{2}$%
\end{tabular}
\begin{tabular}{r}
$n_{1}$%
\end{tabular}
\begin{tabular}{r}
$n_{2}$%
\end{tabular}
\  & (1) & (2) & (3) & (5) & (6) \\ \hline
\begin{tabular}{c}
0.1 \\
\\
\\
\\
\\
0.5 \\
\\
\\
\\
\\
1 \\
\\
\\
\\
\\
2.5
\end{tabular}
\begin{tabular}{r}
5 \\
25 \\
30 \\
50 \\
\\
5 \\
25 \\
30 \\
50 \\
\\
5 \\
25 \\
30 \\
50 \\
\\
5 \\
25 \\
30 \\
50
\end{tabular}
\begin{tabular}{r}
10 \\
25 \\
35 \\
50 \\
\\
10 \\
25 \\
35 \\
50 \\
\\
10 \\
25 \\
35 \\
50 \\
\\
10 \\
25 \\
35 \\
50
\end{tabular}
&
\begin{tabular}{l}
0.071 \\
0.075 \\
0.051 \\
0.046 \\
\\
0.065 \\
0.083 \\
0.056 \\
0.048 \\
\\
0.082 \\
0.075 \\
0.054 \\
0.055 \\
\\
0.092 \\
0.061 \\
0.068 \\
0.048
\end{tabular}
&
\begin{tabular}{l}
0.233 \\
0.116 \\
0.081 \\
0.067 \\
\\
0.274 \\
0.147 \\
0.122 \\
0.095 \\
\\
0.331 \\
0.178 \\
0.148 \\
0.113 \\
\\
0.397 \\
0.208 \\
0.177 \\
0.124
\end{tabular}
&
\begin{tabular}{l}
0.099 \\
0.086 \\
0.059 \\
0.052 \\
\\
0.106 \\
0.096 \\
0.069 \\
0.059 \\
\\
0.141 \\
0.092 \\
0.062 \\
0.061 \\
\\
0.179 \\
0.085 \\
0.078 \\
0.057
\end{tabular}
&
\begin{tabular}{l}
0.035 \\
0.059 \\
0.046 \\
0.043 \\
\\
0.042 \\
0.054 \\
0.054 \\
0.041 \\
\\
0.036 \\
0.051 \\
0.046 \\
0.044 \\
\\
0.034 \\
0.047 \\
0.051 \\
0.047
\end{tabular}
&
\begin{tabular}{l}
0.055 \\
0.071 \\
0.055 \\
0.045 \\
\\
0.051 \\
0.069 \\
0.051 \\
0.045 \\
\\
0.054 \\
0.066 \\
0.046 \\
0.045 \\
\\
0.063 \\
0.059 \\
0.064 \\
0.049
\end{tabular}
\\ \hline
\end{tabular}
\end{center}
\end{table}

\begin{table}
\begin{center}
\caption{ Simulated powers of the tests for $H_{\circ }:\varphi
=1$\ $vs$ $H_{1}:\varphi \neq 1$ at 5\% significance level when
$\mu =0.2$ and $\sigma _{1}^{2}=1.$} \label{table2.5}
\begin{tabular}{|c|ccccc|}
\hline
\begin{tabular}{r}
$\sigma _{2}^{2}$%
\end{tabular}
\begin{tabular}{r}
$n_{1}$%
\end{tabular}
\
\begin{tabular}{r}
$n_{2}$%
\end{tabular}
\  & (1) & (2) & (3) & (5) & (6) \\ \hline
\begin{tabular}{c}
0.1 \\
\\
\\
\\
\\
0.5 \\
\\
\\
\\
\\
1 \\
\\
\\
\\
\\
2.5
\end{tabular}
\begin{tabular}{r}
5 \\
25 \\
30 \\
50 \\
\\
5 \\
25 \\
30 \\
50 \\
\\
5 \\
25 \\
30 \\
50 \\
\\
5 \\
25 \\
30 \\
50
\end{tabular}
\begin{tabular}{r}
10 \\
25 \\
35 \\
50 \\
\\
10 \\
25 \\
35 \\
50 \\
\\
10 \\
25 \\
35 \\
50 \\
\\
10 \\
25 \\
35 \\
50
\end{tabular}
&
\begin{tabular}{l}
0.528 \\
0.909 \\
0.964 \\
0.995 \\
\\
0.171 \\
0.381 \\
0.464 \\
0.631 \\
\\
0.124 \\
0.225 \\
0.280 \\
0.423 \\
\\
0.108 \\
0.199 \\
0.215 \\
0.306
\end{tabular}
&
\begin{tabular}{l}
0.396 \\
0.831 \\
0.933 \\
0.989 \\
\\
0.158 \\
0.157 \\
0.215 \\
0.395 \\
\\
0.190 \\
0.063 \\
0.087 \\
0.188 \\
\\
0.219 \\
0.073 \\
0.048 \\
0.101
\end{tabular}
&
\begin{tabular}{l}
0.539 \\
0.907 \\
0.961 \\
0.955 \\
\\
0.156 \\
0.327 \\
0.403 \\
0.585 \\
\\
0.128 \\
0.199 \\
0.229 \\
0.376 \\
\\
0.124 \\
0.155 \\
0.148 \\
0.247
\end{tabular}
&
\begin{tabular}{l}
0.447 \\
0.891 \\
0.956 \\
0.995 \\
\\
0.156 \\
0.385 \\
0.458 \\
0.633 \\
\\
0.107 \\
0.225 \\
0.280 \\
0.428 \\
\\
0.068 \\
0.193 \\
0.219 \\
0.302
\end{tabular}
&
\begin{tabular}{l}
0.435 \\
0.882 \\
0.952 \\
0.995 \\
\\
0.148 \\
0.365 \\
0.439 \\
0.608 \\
\\
0.109 \\
0.239 \\
0.267 \\
0.417 \\
\\
0.076 \\
0.189 \\
0.201 \\
0.283
\end{tabular}
\\ \hline
\end{tabular}

\end{center}
\end{table}

\begin{table}
\begin{center}
\caption{ Simulated powers of the tests for $H_{\circ }:\varphi
=1$\ $vs$ $H_{1}:\varphi \neq 1$ at 5\% significance level when
$\mu =1$ and $\sigma _{1}^{2}=1.$} \label{table2.6}
\begin{tabular}{|c|ccccc|}
\hline
\begin{tabular}{r}
$\sigma _{2}^{2}$%
\end{tabular}
\begin{tabular}{r}
$n_{1}$%
\end{tabular}
\
\begin{tabular}{r}
$n_{2}$%
\end{tabular}
\  & (1) & (2) & (3) & (5) & (6) \\ \hline
\begin{tabular}{c}
0.1 \\
\\
\\
\\
\\
0.5 \\
\\
\\
\\
\\
1 \\
\\
\\
\\
\\
2.5
\end{tabular}
\begin{tabular}{r}
5 \\
25 \\
30 \\
50 \\
\\
5 \\
25 \\
30 \\
50 \\
\\
5 \\
25 \\
30 \\
50 \\
\\
5 \\
25 \\
30 \\
50
\end{tabular}
\begin{tabular}{r}
10 \\
25 \\
35 \\
50 \\
\\
10 \\
25 \\
35 \\
50 \\
\\
10 \\
25 \\
35 \\
50 \\
\\
10 \\
25 \\
35 \\
50
\end{tabular}
&
\begin{tabular}{l}
1.000 \\
1.000 \\
1.000 \\
1.000 \\
\\
0.999 \\
1.000 \\
1.000 \\
1.000 \\
\\
0.958 \\
1.000 \\
1.000 \\
1.000 \\
\\
0.749 \\
1.000 \\
1.000 \\
1.000
\end{tabular}
&
\begin{tabular}{l}
1.000 \\
1.000 \\
1.000 \\
1.000 \\
\\
0.861 \\
1.000 \\
1.000 \\
1.000 \\
\\
0.559 \\
1.000 \\
1.000 \\
1.000 \\
\\
0.186 \\
0.924 \\
0.971 \\
0.997
\end{tabular}
&
\begin{tabular}{l}
1.000 \\
1.000 \\
1.000 \\
1.000 \\
\\
0.999 \\
1.000 \\
1.000 \\
1.000 \\
\\
0.946 \\
1.000 \\
1.000 \\
1.000 \\
\\
0.702 \\
1.000 \\
1.000 \\
1.000
\end{tabular}
&
\begin{tabular}{l}
1.000 \\
1.000 \\
1.000 \\
1.000 \\
\\
0.981 \\
1.000 \\
1.000 \\
1.000 \\
\\
0.922 \\
1.000 \\
1.000 \\
1.000 \\
\\
0.691 \\
0.998 \\
1.000 \\
1.000
\end{tabular}
&
\begin{tabular}{l}
1.000 \\
1.000 \\
1.000 \\
1.000 \\
\\
0.983 \\
1.000 \\
1.000 \\
1.000 \\
\\
0.927 \\
1.000 \\
1.000 \\
1.000 \\
\\
0.686 \\
0.998 \\
1.000 \\
1.000
\end{tabular}
\\ \hline
\end{tabular}
\end{center}
\end{table}

\begin{table}
\begin{center}
\caption{ Simulated coverage probabilities at 5\% significance
level when $\mu =0$ and $\sigma _{1}^{2}=1.$} \label{table2.7}
\begin{tabular}{|c|ccccc|}
\hline
\begin{tabular}{r}
$\sigma _{2}^{2}$%
\end{tabular}
\begin{tabular}{r}
$n_{1}$%
\end{tabular}
\
\begin{tabular}{r}
$n_{2}$%
\end{tabular}
\  & (2) & (3) & (4) & (5) & (6) \\ \hline
\begin{tabular}{c}
0.1 \\
\\
\\
\\
\\
0.5 \\
\\
\\
\\
\\
1 \\
\\
\\
\\
\\
2.5
\end{tabular}
\begin{tabular}{r}
5 \\
25 \\
30 \\
50 \\
\\
5 \\
25 \\
30 \\
50 \\
\\
5 \\
25 \\
30 \\
50 \\
\\
5 \\
25 \\
30 \\
50
\end{tabular}
\begin{tabular}{r}
10 \\
25 \\
35 \\
50 \\
\\
10 \\
25 \\
35 \\
50 \\
\\
10 \\
25 \\
35 \\
50 \\
\\
10 \\
25 \\
35 \\
50
\end{tabular}
&
\begin{tabular}{l}
0.774 \\
0.884 \\
0.919 \\
0.933 \\
\\
0.726 \\
0.853 \\
0.878 \\
0.905 \\
\\
0.669 \\
0.822 \\
0.852 \\
0.887 \\
\\
0.603 \\
0.792 \\
0.823 \\
0.876
\end{tabular}
&
\begin{tabular}{l}
0.901 \\
0.914 \\
0.941 \\
0.952 \\
\\
0.894 \\
0.904 \\
0.931 \\
0.942 \\
\\
0.859 \\
0.908 \\
0.938 \\
0.942 \\
\\
0.821 \\
0.915 \\
0.922 \\
0.943
\end{tabular}
&
\begin{tabular}{l}
0.743 \\
0.874 \\
0.897 \\
0.914 \\
\\
0.735 \\
0.865 \\
0.884 \\
0.907 \\
\\
0.703 \\
0.856 \\
0.874 \\
0.903 \\
\\
0.642 \\
0.813 \\
0.842 \\
0.882
\end{tabular}
&
\begin{tabular}{l}
0.963 \\
0.939 \\
0.952 \\
0.956 \\
\\
0.957 \\
0.945 \\
0.946 \\
0.959 \\
\\
0.964 \\
0.949 \\
0.953 \\
0.955 \\
\\
0.962 \\
0.953 \\
0.947 \\
0.952
\end{tabular}
&
\begin{tabular}{l}
0.944 \\
0.929 \\
0.945 \\
0.955 \\
\\
0.947 \\
0.954 \\
0.939 \\
0.955 \\
\\
0.943 \\
0.935 \\
0.954 \\
0.954 \\
\\
0.937 \\
0.943 \\
0.938 \\
0.949
\end{tabular}
\\ \hline
\end{tabular}
\end{center}
\end{table}

\begin{table}
\begin{center}
\caption{ Simulated coverage probabilities at 5\% significance
level when $\mu =0.2$ and $\sigma _{1}^{2}=1.$} \label{table2.8}
\begin{tabular}{|c|ccccc|}
\hline
\begin{tabular}{r}
$\sigma _{2}^{2}$%
\end{tabular}
\begin{tabular}{r}
$n_{1}$%
\end{tabular}
\
\begin{tabular}{r}
$n_{2}$%
\end{tabular}
\  & (2) & (3) & (4) & (5) & (6) \\ \hline
\begin{tabular}{c}
0.1 \\
\\
\\
\\
\\
0.5 \\
\\
\\
\\
\\
1 \\
\\
\\
\\
\\
2.5
\end{tabular}
\begin{tabular}{r}
5 \\
25 \\
30 \\
50 \\
\\
5 \\
25 \\
30 \\
50 \\
\\
5 \\
25 \\
30 \\
50 \\
\\
5 \\
25 \\
30 \\
50
\end{tabular}
\begin{tabular}{r}
10 \\
25 \\
35 \\
50 \\
\\
10 \\
25 \\
35 \\
50 \\
\\
10 \\
25 \\
35 \\
50 \\
\\
10 \\
25 \\
35 \\
50
\end{tabular}
&
\begin{tabular}{l}
0.778 \\
0.884 \\
0.924 \\
0.933 \\
\\
0.686 \\
0.853 \\
0.878 \\
0.905 \\
\\
0.661 \\
0.840 \\
0.857 \\
0.888 \\
\\
0.644 \\
0.814 \\
0.831 \\
0.873
\end{tabular}
&
\begin{tabular}{l}
0.901 \\
0.914 \\
0.941 \\
0.951 \\
\\
0.872 \\
0.904 \\
0.931 \\
0.941 \\
\\
0.852 \\
0.931 \\
0.928 \\
0.938 \\
\\
0.851 \\
0.924 \\
0.935 \\
0.936
\end{tabular}
&
\begin{tabular}{l}
0.748 \\
0.877 \\
0.904 \\
0.914 \\
\\
0.724 \\
0.865 \\
0.884 \\
0.907 \\
\\
0.692 \\
0.855 \\
0.882 \\
0.915 \\
\\
0.892 \\
0.841 \\
0.855 \\
0.887
\end{tabular}
&
\begin{tabular}{l}
0.963 \\
0.939 \\
0.947 \\
0.959 \\
\\
0.969 \\
0.945 \\
0.946 \\
0.959 \\
\\
0.965 \\
0.952 \\
0.959 \\
0.944 \\
\\
0.963 \\
0.942 \\
0.946 \\
0.943
\end{tabular}
&
\begin{tabular}{l}
0.944 \\
0.929 \\
0.945 \\
0.955 \\
\\
0.936 \\
0.929 \\
0.949 \\
0.954 \\
\\
0.935 \\
0.946 \\
0.948 \\
0.941 \\
\\
0.937 \\
0.938 \\
0.946 \\
0.941
\end{tabular}
\\ \hline
\end{tabular}
\end{center}
\end{table}

\begin{table}
\begin{center}
\caption{ Simulated coverage probabilities at 5\% significance
level when $\mu =1$ and $\sigma _{1}^{2}=1.$} \label{table2.9}
\begin{tabular}{|c|ccccc|}
\hline
\begin{tabular}{r}
$\sigma _{2}^{2}$
\end{tabular}
\begin{tabular}{r}
$n_{1}$%
\end{tabular}
\
\begin{tabular}{r}
$n_{2}$%
\end{tabular}
\  & (2) & (3) & (4) & (5) & (6) \\ \hline
\begin{tabular}{c}
0.1 \\
\\
\\
\\
\\
0.5 \\
\\
\\
\\
\\
1 \\
\\
\\
\\
\\
2.5
\end{tabular}
\begin{tabular}{r}
5 \\
25 \\
30 \\
50 \\
\\
5 \\
25 \\
30 \\
50 \\
\\
5 \\
25 \\
30 \\
50 \\
\\
5 \\
25 \\
30 \\
50
\end{tabular}
\begin{tabular}{r}
10 \\
25 \\
35 \\
50 \\
\\
10 \\
25 \\
35 \\
50 \\
\\
10 \\
25 \\
35 \\
50 \\
\\
10 \\
25 \\
35 \\
50
\end{tabular}
&
\begin{tabular}{l}
0.771 \\
0.884 \\
0.924 \\
0.929 \\
\\
0.699 \\
0.853 \\
0.884 \\
0.896 \\
\\
0.667 \\
0.827 \\
0.870 \\
0.884 \\
\\
0.614 \\
0.722 \\
0.821 \\
0.876
\end{tabular}
&
\begin{tabular}{l}
0.901 \\
0.914 \\
0.941 \\
0.942 \\
\\
0.868 \\
0.904 \\
0.931 \\
0.938 \\
\\
0.839 \\
0.917 \\
0.932 \\
0.937 \\
\\
0.841 \\
0.915 \\
0.930 \\
0.943
\end{tabular}
&
\begin{tabular}{l}
0.743 \\
0.877 \\
0.904 \\
0.912 \\
\\
0.728 \\
0.865 \\
0.897 \\
0.908 \\
\\
0.724 \\
0.855 \\
0.889 \\
0.894 \\
\\
0.658 \\
0.813 \\
0.854 \\
0.882
\end{tabular}
&
\begin{tabular}{l}
0.963 \\
0.939 \\
0.947 \\
0.945 \\
\\
0.958 \\
0.945 \\
0.951 \\
0.946 \\
\\
0.958 \\
0.948 \\
0.962 \\
0.949 \\
\\
0.961 \\
0.953 \\
0.956 \\
0.952
\end{tabular}
&
\begin{tabular}{l}
0.944 \\
0.929 \\
0.949 \\
0.947 \\
\\
0.936 \\
0.929 \\
0.943 \\
0.944 \\
\\
0.937 \\
0.939 \\
0.951 \\
0.951 \\
\\
0.932 \\
0.941 \\
0.951 \\
0.949
\end{tabular}
\\ \hline
\end{tabular}
\end{center}
\end{table}

\end{document}